# A simple adaptive estimator of the integrated square of a density

EVARIST GINÉ[*] and RICHARD NICKL[**]

*Department of Mathematics, University of Connecticut, Storrs, CT 06269-3009, USA.
E-mail: [*]gine@math.uconn.edu; [**]nickl@math.uconn.edu*

Given an i.i.d. sample $X_1, \ldots, X_n$ with common bounded density $f_0$ belonging to a Sobolev space of order $\alpha$ over the real line, estimation of the quadratic functional $\int_{\mathbb{R}} f_0^2(x)\,\mathrm{d}x$ is considered. It is shown that the simplest kernel-based plug-in estimator

$$\frac{2}{n(n-1)h_n} \sum_{1 \le i < j \le n} K\Big(\frac{X_i - X_j}{h_n}\Big)$$

is asymptotically efficient if $\alpha > 1/4$ and rate-optimal if $\alpha \le 1/4$. A data-driven rule to choose the bandwidth $h_n$ is then proposed, which does not depend on prior knowledge of $\alpha$, so that the corresponding estimator is rate-adaptive for $\alpha \le 1/4$ and asymptotically efficient if $\alpha > 1/4$.

## 1. Introduction

The estimation of a quadratic functional of a density $f_0$, in particular of $\int f_0^2$, has attracted much interest in the literature since Bickel and Ritov (1988) showed that such functionals can be estimated at the rate $1/\sqrt{n}$ if $f_0$ is $\alpha$-Hölder of order $\alpha > 1/4$ and that this rate cannot be achieved if $\alpha < 1/4$. Such functionals have several statistical applications. For instance, $\int f_0^2$ occurs in Taylor expansions of more complex integral functionals, such as the entropy $\int f_0 \log f_0$; see, for example, Laurent (1996). They are also part of constants appearing in the exact expression of the MISE of kernel density estimators and hence their estimates can be used in optimal bandwidth selection. Bickel and Ritov (1988) constructed an efficient and $\sqrt{n}$-consistent kernel-based estimator for $\int f_0^2$ and Laurent (1996) achieved the same for an estimator based on orthogonal series. The treatment of the bias term by these authors necessitated rather complicated expressions for the actual estimators, which consist of the difference of two U-statistics. As a first goal of this article, we show that the simplest 'plug-in' kernel density estimator introduced in Hall and Marron (1987),

$$T_n(h_n) := \frac{2}{n(n-1)h_n} \sum_{1 \leq i < j \leq n} K\left(\frac{X_i - X_j}{h_n}\right), \tag{1}$$

where $X_i$ are i.i.d. with common density $f_0$ on the real line, also does the job. ($T_n$ is obtained as follows. Estimate $f_0$ by the usual kernel density estimator and estimate integration with respect to $f_0$ by integration with respect to the empirical measure, then delete the diagonal terms.) Our main point here consists of an observation on the bias term based on smoothing properties of convolutions, borrowed in part from Giné and Nickl (2007). We note that in the context-of-goodness of fit tests, Butucea (2007) considered a different (but related) kernel-based estimator for $\int f_0^2$, where $K(x)$ must equal $\sin(x)/\pi x$. Our results also hold for her estimator, without any (other than the usual) restrictions on the kernel; see Remark 1 below. The same methods can be applied, with the natural changes, to other quadratic functionals, such as $\int (f_0^{(k)})^2$ for $k > 0$.

As is well known, efficient estimation of $\int f_0^2$ is possible if $f_0$ is in a Sobolev space of order $\alpha > 1/4$, but in the 'low regularity case' $\alpha \leq 1/4$, the best rate of convergence is $n^{-4\alpha/(4\alpha+1)}$. We show that $T_n(h_n)$ achieves this rate if one chooses the bandwidth $h_n$ of the right order, where $h_n$ depends on the unknown quantity $\alpha$. It is then natural to ask whether one can choose the bandwidth in some data-dependent way, so as to obtain an estimator of $\int f_0^2$ which is rate-adaptive over Sobolev balls if $\alpha \leq 1/4$ and efficient if $\alpha > 1/4$. Using Lepski's method (Lepski (1990), Lepski and Spokoiny (1997)), we show that this is in fact possible for the simple estimator $T_n(h_n)$. Rate-adaptive estimation of $\int f_0^2$ was first considered by Efromovich and Low (1996), and more recently by, for example, Laurent and Massart (2000), Laurent (2005), Cai and Low (2006) and Klemelä (2006). None of these authors used kernel-based estimators, and, except for Laurent (2005), all of them worked in the context of a Gaussian white noise model. Since we are interested in the low-regularity case where $\alpha < 1/4$, the restriction to the Gaussian white noise model is inconvenient, as it is not clear how asymptotic results in the Gaussian white noise model translate into the usual density model in this case. It turns out that deriving our results in the more general density model on the real line leads to no major complications. Our derivations rely on elementary U-statistic theory, some simple Fourier analytical methods and a recent exponential inequality for canonical U-statistics of order 2 due to Giné, Latała and Zinn (2000), with constants obtained in Houdré and Reynaud-Bouret (2003). A discussion of the relationship of our results to those in Laurent (2005) is given in Remark 5 below.

## 2. Basic setup

We will assume that the probability density $f_0$ is bounded, that is, $f_0 \in L^\infty := L^\infty(\mathbb{R})$, and contained in a Sobolev space of order $\alpha > 0$, defined as follows. First, denote by $L^p := L^p(\mathbb{R}, \lambda)$ the usual spaces of measurable functions $\phi$ satisfying $\|\phi\|_p^p := \int_\mathbb{R} |\phi(x)|^p \, dx < \infty$ for $1 \leq p < \infty$. For $\phi \in L^1$, we define the Fourier-transform by $F\phi(u) = \int_\mathbb{R} e^{-ixu} \phi(x) \, dx$ and extend it continuously to $L^2$. ($F$ is, up to a multiplicative constant, the Fourier–Plancherel transform.) We then set

$$H_2^\alpha = H_2^\alpha(\mathbb{R}) := \{\phi \in L^2 : \|\phi\|_{2,\alpha} = \|F\phi(\cdot)(1+|\cdot|^2)^{\alpha/2}\|_2 < \infty\}.$$

We note that a common equivalent characterization of $H_2^\alpha$ is in terms of integrated $L^2$-Hölder conditions: for $\phi \in L^2$ and $0 < \alpha < 1$, define

$$I_\alpha(\phi) = \int_\mathbb{R} \int_\mathbb{R} \frac{|\phi(x-t) - \phi(x)|^2}{|t|^{1+2\alpha}} \, dx \, dt.$$

It can then be shown that $\phi \in H_2^\alpha$ if and only if $\phi \in L^2$ and $I_\alpha(\phi) < \infty$ (cf. page 144 in Malliavin (1995)). Throughout the proofs, we will freely use these and other standard facts from Fourier analysis, as well as Young's inequalities for convolutions. We refer, for example, to Chapter III in Malliavin (1995) or Chapter 8 in Folland (1999). Also, unless otherwise indicated, all integrals in this article will be over the real line.

It is also convenient to introduce U-statistic notation. For a symmetric function of two variables $R(x,y)$, we write

$$U_n^{(2)}(R) = \frac{2}{n(n-1)} \sum_{1 \leq i < j \leq n} R(X_i, X_j).$$

We recall (e.g., de la Peña and Giné (1999), page 137) that the Hoeffding projections of $R$ are

$$\pi_1 R(x) = ER(x, X_1) - ER(X_1, X_2),$$
$$\pi_2(R)(x,y) = R(x,y) - ER(x, X_1) - ER(y, X_1) + ER(X_1, X_2),$$

which induce the Hoeffding decomposition

$$U_n^{(2)}(R) - ER(X_1, X_2) = 2U_n^{(1)}(\pi_1 R) + U_n^{(2)}(\pi_2 R), \tag{2}$$

where $U_n^{(1)}(\pi_1 R) = n^{-1} \sum_{i=1}^n (\pi_1 R)(X_i)$. Note that, by orthogonality,

$$E(U_n^{(1)}(\pi_1 R))^2 = n^{-1} E((\pi_1 R)(X_1))^2,$$
$$E(U_n^{(2)}(\pi_2 R))^2 = \frac{2}{n(n-1)} E((\pi_2 R)(X_1, X_2))^2.$$

# 3. Estimation of $\int_\mathbb{R} f_0^2(x) \, dx$

The simple estimator $T_n(h_n)$ defined in (1) can be shown to be optimal, as we prove in this section.

Here and elsewhere in this article, we take the kernel $K$ in (1) to be a symmetric and bounded function such that $\int K(u) \, du = 1$, as well as $\int |K(u)||u| \, du < \infty$ and $0 < h_n \to 0$.

For ease of notation, we will often write $K_{h_n}(x) := h_n^{-1} K(x/h_n)$. Also, we define the Sobolev ball $\mathcal{H}_\alpha(R) = \{\phi : \|\phi\|_{2,\alpha} \leq R\}$ and $\mathcal{B}(L) = \{\phi : \|\phi\|_\infty \leq L\}$.

**Theorem 1.** *Let $f_0 \in H_2^\alpha \cap L^\infty$ for some $0 < \alpha \leq 1/2$.*
I. *We have*

$$\sup_{f_0 \in \mathcal{H}_\alpha(R)} \left| ET_n(h_n) - \int_{\mathbb{R}} f_0^2(x)\,dx \right| \leq B(h_n) := c_1(R) h_n^{2\alpha} \tag{3}$$

*and*

$$\sup_{f_0 \in \mathcal{H}_\alpha(R) \cap \mathcal{B}(L)} E\left( T_n(h_n) - ET_n(h_n) - \frac{1}{n}\sum_{i=1}^{n} Y_i \right)^2 \leq c_2^2(R) \sigma^2(h_n, n)$$
$$:= c_2^2(R) \left( \frac{1}{n^2 h_n} \vee \frac{L h_n^{2\alpha}}{n} \right), \tag{4}$$

*where $Y_i = 2(f_0(X_i) - \int_{\mathbb{R}} f_0^2)$ and where $c_1(R)$ and $c_2(R)$ are numerical constants depending only on $R$ and the function $K$.*
II. *As a consequence, taking $h_n$ so that $h_n \approx n^{-2/(4\alpha+1)}$, we have the following:*

(a) *if $0 < \alpha \leq 1/4$, then*

$$T_n(h_n) - \int_{\mathbb{R}} f_0(x)^2\,dx = O_P(n^{-4\alpha/(4\alpha+1)});$$

(b) *if $\alpha > 1/4$, and if $\tau^2 = [\int_{\mathbb{R}} f_0^3 - (\int_{\mathbb{R}} f_0^2)^2]$, then*

$$\sqrt{n}\left( T_n(h_n) - \int_{\mathbb{R}} f_0(x)^2\,dx \right) \to_d Z \sim N(0, 4\tau^2).$$

**Proof.** We first treat the bias term, where we adapt an observation due to Giné and Nickl (2007), Section 4.1.1, to the present situation. The bias equals

$$ET_n(h_n) - \int f_0^2 = \int_{\mathbb{R}} \int_{\mathbb{R}} K_{h_n}(x-y) f_0(y)\,dy\, f_0(x)\,dx - \int_{\mathbb{R}} f_0(x) f_0(x)\,dx$$
$$= \int_{\mathbb{R}} \int_{\mathbb{R}} K_{h_n}(x-y)[f_0(y) - f_0(x)] f_0(x)\,dy\,dx$$
$$= \int_{\mathbb{R}} \int_{\mathbb{R}} K(u)[f_0(x - uh_n) - f_0(x)] f_0(x)\,du\,dx \tag{5}$$
$$= \int_{\mathbb{R}} K(u) \left[ \int_{\mathbb{R}} \bar{f}_0(uh_n - x) f_0(x)\,dx - \int_{\mathbb{R}} \bar{f}_0(0 - x) f_0(x)\,dx \right] du$$
$$= \int_{\mathbb{R}} K(u)[(\bar{f}_0 * f_0)(uh_n) - (\bar{f}_0 * f_0)(0)]\,du,$$

where $\bar{f}_0(x) = f_0(-x)$ and $*$ denotes convolution. The essential observation now is that the smoothness of $\bar{f}_0 * f_0$ will be of order $2\alpha$ instead of just $\alpha$, due to the smoothing properties of convolutions. The following elementary Fourier analytic lemma shows how this applies in our setup.

**Lemma 1.** *Suppose that $f, g \in H_2^\alpha$ with $0 < \alpha \le 1/2$. Then, for any $x \in \mathbb{R}$ and $t \ne 0$,*

$$\frac{|(f*g)(x+t) - (f*g)(x)|}{|t|^{2\alpha}} \le C \|f\|_{2,\alpha} \|g\|_{2,\alpha},$$

*where $0 < C < \infty$ is a fixed constant that does not depend on $f, g, x$ or $t$.*

**Proof.** As we will only use this lemma for $f, g \in L^1$, and in order to avoid some technicalities, we will prove it only in this case. Hence, $f * g$ is in $L^1$ and is continuous and, since $f, g$ are also in $L^2$, we also have $F(f * g) \in L^1$. Consequently, we can apply the Fourier inversion theorem to obtain

$$\begin{aligned}
\frac{|(f*g)(x+t) - (f*g)(x)|}{|t|^{2\alpha}} &\le |t|^{-2\alpha} \|F^{-1} F[(f*g)(\cdot+t) - (f*g)(\cdot)]\|_\infty \\
&\le (2\pi)^{-1} |t|^{-2\alpha} \|F[(f*g)(\cdot+t) - (f*g)(\cdot)]\|_1 \\
&= (2\pi)^{-1} |t|^{-2\alpha} \int_{\mathbb{R}} |F(f*g)(u)[e^{-iut} - 1]| \, du \\
&= (2\pi)^{-1} \int_{\mathbb{R}} |Ff| |u|^\alpha |Fg| |u|^\alpha \frac{|e^{-iut} - e^{-i0}|}{|u|^{2\alpha} |t|^{2\alpha}} \, du \\
&\le C \|f\|_{2,\alpha} \|g\|_{2,\alpha}
\end{aligned}$$

since $e^{-i(\cdot)}$ is bounded Lipschitz. $\square$

This lemma and identity (5) now give, by the conditions on the kernel, that

$$\left| ET_n(h_n) - \int f_0^2 \right| \le c_1 h_n^{2\alpha},$$

where $c_1 = C\|f_0\|_{2,\alpha}^2 \int |K(u)| |u|^{2\alpha} \, du \le CR^2 \int |K(u)| |u|^{2\alpha} \, du$, that is, (3).

Next, we show (4). Setting

$$R(u,v) := K_{h_n}(u - v),$$

we can write, in U-statistic notation, $T_n(h_n) = U_n^{(2)}(R)$ or, if $\tilde{R}(u,v) = R(u,v) - ER(X_1, X_2)$,

$$T_n(h_n) - ET_n(h_n) = U_n^{(2)}(\tilde{R}).$$

So, by Hoeffding's decomposition (2), it remains to estimate the following statistics (note that $\pi_i R = \pi_i \tilde{R}$, $i=1,2$):

$$U_n^{(2)}(\tilde{R}) - \frac{1}{n}\sum_{i=1}^n Y_i = \left(2U_n^{(1)}(\pi_1 R) - \frac{1}{n}\sum_{i=1}^n Y_i\right) + U_n^{(2)}(\pi_2 R) =: S_1 + S_2.$$

First, we have, by Plancherel,

$$\begin{aligned}
nES_1^2 &\leq E\left[\int 2K_{h_n}(X_1 - y)f_0(y)\,\mathrm{d}y - 2f_0(X_1)\right]^2 \\
&\leq 4\|f_0\|_\infty \|K_{h_n} * f_0 - f_0\|_2^2 \\
&= \frac{4}{\sqrt{2\pi}}\|f_0\|_\infty \|(FK_{h_n} - 1)\cdot|\cdot|^{-\alpha}Ff_0\cdot|\cdot|^\alpha\|_2^2 \\
&\leq \frac{4}{\sqrt{2\pi}}\|f_0\|_\infty h_n^{2\alpha}\left(\sup_u \frac{|FK(h_n u) - FK(0)|}{|uh_n|^\alpha}\right)^2 \|f_0\|_{2,\alpha}^2 \\
&= \frac{4}{\sqrt{2\pi}}\|f_0\|_\infty h_n^{2\alpha}\sup_u\left|\int \frac{\mathrm{e}^{-\mathrm{i}xh_n u} - 1}{|uh_n|^\alpha}K(x)\,\mathrm{d}x\right|^2 \|f_0\|_{2,\alpha}^2 \\
&\leq \frac{8}{\sqrt{2\pi}}\|f_0\|_\infty \left(\int |K(x)||x|^\alpha\,\mathrm{d}x\right)^2 \|f_0\|_{2,\alpha}^2 h_n^{2\alpha}.
\end{aligned} \tag{6}$$

Next, since $\pi_2$ is a projection of $L^2(f_0(x)\,\mathrm{d}x)$, it follows from Young's inequalities that

$$\begin{aligned}
ES_2^2 &\leq \frac{2}{n(n-1)}ER^2 = \frac{2}{n(n-1)}E[K_{h_n}(X_1 - X_2)]^2 \\
&= \frac{2}{n(n-1)}\int (K_{h_n}^2 * f_0)(y)f_0(y)\,\mathrm{d}y \\
&\leq \frac{2\|f_0\|_2^2 \|K\|_2^2}{n(n-1)h_n}.
\end{aligned} \tag{7}$$

Now, (6) and (7) complete the proof of (4). The remaining claims in Part II follow by the choice of the bandwidth and, in case (a) (and hence $\alpha \leq 1/4$), noting that we have $n^{-1}\sum_{i=1}^n Y_i = O_P(n^{-1/2}) = O_P(n^{-4\alpha/(4\alpha+1)})$ and, in case (b), from the central limit theorem for the random variables $Y_i$. $\square$

Without loss of generality, we restricted ourselves to $0 < \alpha \leq 1/2$ in Theorem 1. It is obvious that Part II holds for all $\alpha > 0$ and it can be seen that Part I does too, although this is not of interest here.

**Remark 1.** A second plug-in estimator of $\int f_0^2$ is

$$\overline{T}_n(h_n) = \frac{2}{n(n-1)} \sum_{1 \leq i < j \leq n} \int K_{h_n}(x - X_i) K_{h_n}(x - X_j) \, dx,$$

obtained by integrating the square of the usual kernel density estimator and deleting the diagonal terms. Although Theorem 1 could also be proved for this estimator, we choose to work with $T_n(h_n)$ because it is simpler to compute. The results of Theorem 1 for $\overline{T}_n$ can be derived by similar computations. Here, we briefly consider the bias, which is really the main part, by relating it to the bias of $T_n$ (as in Bickel and Ritov (1988)): using (3) and (6), we have

$$\left| E\overline{T}_n(h_n) - \int f_0^2 \right| = \left| E\overline{T}_n(h_n) - 2ET_n(h_n) + \int f_0^2 + 2\left(ET_n(h_n) - \int f_0^2\right) \right|$$

$$\leq \left| E\overline{T}_n(h_n) - 2ET_n(h_n) + \int f_0^2 \right| + 2c_1 h_n^{2\alpha}$$

$$= \int [(K_{h_n} * f_0)(x) - f_0(x)]^2 \, dx + 2c_1 h_n^{2\alpha}$$

$$\leq \frac{2}{\sqrt{2\pi}} \left( \int |K(u)| |u|^\alpha \, du \right)^2 \|f_0\|_{2,\alpha}^2 h_n^{2\alpha} + 2c_1 h_n^{2\alpha}.$$

Butucea (2007) also obtains such a bound, but only for the special kernel $K(x) = \sin(x)/\pi x$, and it is the use of Lemma 1 that allows us to consider the case of general kernels.

**Remark 2.** Bickel and Ritov (1988) show that if $|f_0(x+h) - f_0(x)| \leq g(x)|h|^\alpha$ for some $g \in L^2 \cap L^\infty$, $x \in \mathbb{R}$, $|h| < 1$ and $\alpha > 1/4$, then

$$\sqrt{n}\left(2T_n(h_n) - \overline{T}_n(h_n) - \int f_0^2\right) \to_d Z \sim N(0, 4\tau^2)$$

with $h_n = n^{-2/(4\alpha+1)}$ (actually they consider a 'decoupled' version). Clearly, any such $f_0$ is contained in $H_2^\beta$ for all $\beta < \alpha$ and this implies, by Theorem 1, that the simpler estimator $T_n(h_n)$ satisfies the same central limit theorem. (As a matter of fact, Lemma 1 and hence Theorem 1 also holds for such $f_0$, even without requiring $g \in L^\infty$; see Lemma 12 (and the discussion following it) in Giné and Nickl (2007). However, the proofs there are much more technical, which is why we prefer to work with Sobolev spaces here.)

## 4. Adaptive estimation of $\int_\mathbb{R} f_0^2(x) \, dx$

In Theorem 1, one must know $\alpha$ in order to choose $h_n$ in an optimal way, $h_n$ ranging between $n^{-2}$ and 1. We will now use $T_n(h_n)$ to construct a kernel-based rate-adaptive estimator of $\int f_0^2$ that requires only that $\int f_0^2$ is bounded by a known constant $L$ and that $f_0$ is a bounded function contained in $H_2^\alpha$ for some (unknown) $\alpha > 0$. In practice, one can restrict oneself to $0 < \alpha < 1/4 + \varepsilon$ with $\varepsilon$ positive and arbitrary since the rate of convergence $n^{-1/2}$ in part (b) of Theorem 1 could not be improved if one knew that $\alpha$ were larger. In particular, it suffices to consider bandwidths that are faster than $n^{-1+\delta}$ for some arbitrarily small $\delta$.

In what follows, we borrow in part from methods developed by Lepski and Spokoiny (1997) for kernel-based pointwise adaptive estimation in the Gaussian white noise model. Our situation, however, is substantially different in several respects. For instance, there is a critical breakpoint in convergence rates at $\alpha = 1/4$ and we do not have the convenience of immediate Gaussian tail inequalities.

For any given $n \in \mathbb{N}$, $n > 1$, we define a grid of bandwidths

$$\mathcal{H} := \left\{ h \in \left[ \frac{(\log n)^4}{n^2}, \frac{1}{n^{1-\delta}} \right] : h_0 = \frac{1}{n^{1-\delta}}, \ h_1 = \frac{\log n}{n}, h_2 = \frac{\ell(n)}{n}, h_{k+1} = \frac{h_k}{\rho}, k = 2, 3, \ldots \right\},$$

where $\rho > 1$ and $\ell(n)$ is any function such that $\ell(n) \to 0$ and $\ell(n) \log n \to \infty$ as $n \to \infty$, and $\ell(n) < \log n$ for all $n$. (In particular, $\ell(n)$ can be chosen to tend to zero as slowly as desired.) It is easy to check that the number of elements in this grid is smaller than $3 + (\log n)/(\log \rho) = O(\log n)$ and we shall use this estimate below. Next, we define the function $d(h)$ for all $h \in [n^{-2}(\log n)^4, n^{-1+\delta}]$ as

$$d(h) = \sqrt{2M \log \frac{h_0}{h}} \quad \text{for } h < h_2 \quad \text{and} \quad d(h) = \ell(n)^{-1/2} \quad \text{for } h_0 \geq h \geq h_2,$$

where $M := 12^2 \|K\|_2^2 L$ and where we recall that $L$ is a bound on $\int f_0^2$. We also set $\tilde{\sigma}(h, n) = n^{-1} h^{-1/2}$. The bandwidth estimator is defined as

$$\hat{h}_n = \max\{h \in \mathcal{H} : |T_n(h) - T_n(g)| \leq \tilde{\sigma}(g, n) d(g) \ \forall g < h, g \in \mathcal{H}\}.$$

***Remark 3.*** If $h$ equals the next to last element in the grid $\mathcal{H}$ and $g$ is the last, then $\tilde{\sigma}(g, n) d(g)$ is of the order $(\log n)^{-3/2}$, whereas $|T_n(h) - T_n(g)| = O_P((\log n)^{-2})$, by Theorem 1. Hence, $\hat{h}_n$ exists with probability tending to 1 as $n \to \infty$. In the next theorem, expectations that involve events based on $\hat{h}_n$ should be understood as taken over the event $\{\hat{h}_n \text{ exists}\}$.

***Remark 4.*** In cases (a) and (b) in Theorem 2 below, the rates of convergence obtained are, in fact, slightly slower than those in Theorem 1. This is not surprising, as Efromovich and Low (1996) showed that one must pay exactly these penalties if one wants to estimate $\int f_0^2$ adaptively.

***Remark 5.*** Laurent (2005) considered adaptive estimation of $\int_{\mathbb{R}} f_0^2$ by model selection. Her results are comparable to our Theorem 2 below. (She considers $f_0$ contained in the Besov space $B_{2,\infty}^{\alpha}$ which is slightly more general in view of the imbeddings $B_{2,\infty}^{\alpha+\varepsilon} \subseteq H_2^{\alpha} \subseteq B_{2,\infty}^{\alpha}$ for every $\varepsilon > 0$.) In her Theorem 1, she assumes that an a priori bound for $\|f_0\|_\infty$ is known. Similarly, we have the assumption of a known upper bound $L$ for $\int f_0^2$. In her Theorem 2, Laurent (2005) proposes a remedy for this problem by estimating this upper bound. Similarly, we could estimate the upper bound $L$ by $T_n(h_{\min})$ to achieve the same goal.

**Theorem 2.** *Let $f_0 \in H_2^{\alpha} \cap L^\infty$ for some $\alpha > 0$.*

(a) *If $0 < \alpha < 1/4$, then*

$$T_n(\hat{h}_n) - \int_{\mathbb{R}} f_0^2(x)\,\mathrm{d}x = O_P\left(\left(\frac{\sqrt{\log n}}{n}\right)^{4\alpha/(4\alpha+1)}\right).$$

(b) *If $\alpha = 1/4$, then*

$$T_n(\hat{h}_n) - \int_{\mathbb{R}} f_0^2(x)\,\mathrm{d}x = O_P(n^{-1/2}\ell(n)^{-1}).$$

(c) *If $\alpha > 1/4$ and $\tau^2 = [\int_{\mathbb{R}} f_0^3 - (\int_{\mathbb{R}} f_0^2)^2]$, then*

$$\sqrt{n}\left(T_n(\hat{h}_n) - \int_{\mathbb{R}} f_0^2(x)\,\mathrm{d}x\right) \to_d Z \sim N(0, 4\tau^2).$$

**Proof.** We first observe that $\tilde{\sigma}(h,n) = \sigma(h,n)$ whenever $h \leq h_1$, which will always be the case in this proof. Define $h_f(=h_{f_0})$ as $h_1$ if $\alpha > 1/4$, as $h_2$ if $\alpha = 1/4$ and, otherwise,

$$h_f = \max\{h \in \mathcal{H} : c_1 h^{2\alpha} \leq \tfrac{1}{4}\sigma(h,n)d(h),\ h < h_2\}.$$

It is easily checked that $h_f$ exists and is of the order of $(n/\sqrt{\log n})^{-2/(4\alpha+1)}$ if $\alpha < 1/4$. By construction in case $\alpha < 1/4$ and by straightforward computations in the other two cases, we have, for $n$ large enough,

$$B(h_f) \leq \tfrac{1}{4}\sigma(h_f,n)d(h_f). \tag{8}$$

We estimate the expectation of

$$\left|T_n(\hat{h}_n) - \int_{\mathbb{R}} f_0^2(x)\,\mathrm{d}x - \frac{1}{n}\sum_{i=1}^n Y_i\right|$$

over each of the two events $\{\hat{h}_n \geq h_f\}$ and $\{\hat{h}_n < h_f\}$. In the first case, we have

$$E\left|T_n(\hat{h}_n) - \int_{\mathbb{R}} f_0^2(x)\,\mathrm{d}x - \frac{1}{n}\sum_{i=1}^n Y_i\right|I_{[\hat{h}_n \geq h_f]}$$

$$\leq E\left[|T_n(\hat{h}_n) - T_n(h_f)| + \left|T_n(h_f) - ET_n(h_f) - \frac{1}{n}\sum_{i=1}^{n}Y_i\right| + \left|ET_n(h_f) - \int f_0^2\right|\right] I_{[\hat{h}_n \geq h_f]}$$

$$\leq \sigma(h_f, n)d(h_f) + c_2\sigma(h_f, n) + B(h_f)$$

$$= O(\sigma(h_f, n)d(h_f)),$$

where we use the definition of $\hat{h}_n$, Theorem 1 and (8). In the other case, where $\{\hat{h}_n < h_f\}$, we will rely on the following lemma, which will be proved below.

**Lemma 2.** *Let $h \in \mathcal{H}$ and $h < h_f$. There exists a constant $D < \infty$ so that, for all $n$ large enough, if $h < h_2$, then,*

$$\Pr(\hat{h}_n = h) \leq D(\log n)\exp(-d^2(h)/M)$$

*and if $h = h_2$, then*

$$\Pr(\hat{h}_n = h_2) \leq D[\exp(-d^2(h_2)/M) + (\log n)\exp(-d^2(h_3)/M)].$$

This lemma, Theorem 1, (8), the size of the grid and the definition of $d(h)$ now give, for $\alpha \leq 1/4$ and hence for $h_f \leq h_2$,

$$E\left(\left|T_n(\hat{h}_n) - \int_{\mathbb{R}} f_0^2(x)\,dx - \frac{1}{n}\sum_{i=1}^{n}Y_i\right| I_{[\hat{h}_n < h_f]}\right)$$

$$= \sum_{h \in \mathcal{H}: h < h_f} E\left(\left|T_n(h) - \int_{\mathbb{R}} f_0^2(x)\,dx - \frac{1}{n}\sum_{i=1}^{n}Y_i\right| I_{[\hat{h}_n = h]}\right) \quad (9)$$

$$\leq \sum_{h \in \mathcal{H}: h < h_f} E\left(\left[\left|T_n(h) - ET_n(h) - \frac{1}{n}\sum_{i=1}^{n}Y_i\right| + \left|ET_n(h) - \int f_0^2\right|\right] I_{[\hat{h}_n = h]}\right)$$

$$\leq \sum_{h \in \mathcal{H}: h < h_f} \left(E\left|T_n(h) - ET_n(h) - \frac{1}{n}\sum_{i=1}^{n}Y_i\right|^2\right)^{1/2} (\Pr(\hat{h}_n = h))^{1/2} + B(h_f)$$

$$\leq c_2 D^{1/2}(\log n)^{1/2} n^{-\delta} \sum_{h \in \mathcal{H}: h < h_f} h^{-1/2} h + \frac{1}{4}\sigma(h_f, n)d(h_f)$$

$$\leq D' n^{-\delta}(\log n)^{3/2} h_f^{1/2} + \frac{1}{4}\sigma(h_f, n)d(h_f)$$

$$= Z_n(\alpha) + O(\sigma(h_f, n)d(h_f)),$$

where $D'$ is an absolute constant and where $Z_n(\alpha) = o(n^{-1/2})$ if $\alpha \geq 1/4$ and

$$Z_n(\alpha) = o\left(\left(\frac{\sqrt{\log n}}{n}\right)^{4\alpha/(4\alpha+1)}\right)$$

otherwise (as can easily be seen from the definition of $h_f$). If $\alpha > 1/4$ (and hence $h_f = h_1$), then one must add the term

$$\left(E\left|T_n(h_2) - ET_n(h_2) - \frac{1}{n}\sum_{i=1}^n Y_i\right|^2\right)^{1/2} (\Pr(\hat{h}_n = h_2))^{1/2}$$
$$\leq D^{1/2}(n\ell(n))^{-1/2}[\exp(-d^2(h_2)/M) + (\log n)\exp(-d^2(h_3)/M)]^{1/2}$$
$$= o(n^{-1/2})$$

in the sum over $h < h_f$ in the line before (9), hence yielding the same result.

Summarizing these findings, we conclude that

$$E\left|T_n(\hat{h}_n) - \int_\mathbb{R} f_0^2(x)\,dx - \frac{1}{n}\sum_{i=1}^n Y_i\right| = O(\sigma(h_f,n)d(h_f)) + Z_n(\alpha) + o(n^{-1/2}). \quad (10)$$

By definition of $h_f$, it follows that, if $\alpha > 1/4$,

$$\sigma(h_f,n)d(h_f) \approx n^{-1+(1/2)}(\log n)^{-1/2}(\ell(n))^{-1/2} = o(n^{-1/2}),$$

hence giving the central limit theorem in part (c) of the theorem by (10). Similarly, for part (b), if $\alpha = 1/4$ and hence $h_f = h_2$, we obtain

$$\sigma(h_f,n)d(h_f) \approx n^{-1+(1/2)}(\ell(n))^{-1/2}(\ell(n))^{-1/2} = O(n^{-1/2}\ell(n)^{-1})$$

and if $\alpha < 1/4$, we have

$$\sigma(h_f,n)d(h_f) \approx \frac{\sqrt{\log n}}{n}\left(\frac{n}{\sqrt{\log n}}\right)^{1/(4\alpha+1)} = O\left(\left(\frac{\sqrt{\log n}}{n}\right)^{4\alpha/(4\alpha+1)}\right),$$

giving part (a). □

It hence remains to prove Lemma 2, where we will use Bernstein's inequality and an exponential inequality for canonical U-statistics of order 2.

**Proof of Lemma 2.** Choose some $h < h_f$, $h \in \mathcal{H}$ and let $h_+ = \rho h$ be the previous element in the grid. By definition of $\hat{h}_n$, we have

$$\Pr(\hat{h}_n = h) \leq \sum_{g \in \mathcal{H}: g \leq h} \Pr(|T_n(g) - T_n(h_+)| > \sigma(g,n)d(g)).$$

However,

$$|T_n(g) - T_n(h_+)| \leq |T_n(g) - ET_n(g) - (T_n(h_+) - ET_n(h_+))| + B(g) + B(h_+),$$

where $g \leq h < h_f$ and also $h_+ \leq h_f$ since $h_f \in \mathcal{H}$. Consequently, by (8),

$$B(g) + B(h_+) \leq 2B(h_f) \leq \tfrac{1}{2}\sigma(h_f,n)d(h_f) \leq \tfrac{1}{2}\sigma(g,n)d(g).$$

Hence,

$$\Pr(\hat{h}_n = h) \leq \sum_{g \in \mathcal{H}\,:\,g \leq h} \Pr(|T_n(g) - ET_n(g) - (T_n(h_+) - ET_n(h_+))| > \tfrac{1}{2}\sigma(g,n)d(g)). \quad (11)$$

For ease of notation, we set

$$L(x,y) := L_g(x,y) = K_g(x-y) - K_{h_+}(x-y)$$

and

$$C_{n,g} := \frac{1}{2}\sigma(g,n)d(g) = \frac{1}{2}\frac{d(g)}{ng^{1/2}}.$$

In particular, in U-statistic notation, we have

$$U_n^{(2)}(L) - EU_n^{(2)}(L) = (T_n(g) - ET_n(g)) - (T_n(h_+) - ET_n(h_+))$$

and, recalling the Hoeffding decomposition (2), we have

$$U_n^{(2)}(L) - EU_n^{(2)}(L) = 2U_n^{(1)}(\pi_1 L) + U_n^{(2)}(\pi_2 L).$$

So, to estimate the right-hand side of (11), it suffices to bound

$$\Pr\{|U_n^{(1)}(\pi_1 L)| > \tau C_{n,g}/2\} \quad \text{and} \quad \Pr\{|U_n^{(2)}(\pi_2 L)| > (1-\tau)C_{n,g}\}$$

for some $0 < \tau < 1$. We will apply Bernstein's inequality (e.g., de la Peña and Giné (1999), page 166) to the linear part (the first probability) and its generalization for canonical U-statistics of order 2 (Giné, Latała and Zinn (2000)) with constants (Houdré and Reynaud-Bouret (2003)) to the second.

*Linear term:* Noting that $\operatorname{Var}(\pi_1 L) \leq E(E_{X_2}L(X_1,X_2))^2$ and that

$$E(E_{X_2}(K_g(X_1 - X_2)))^2 = \int (K_g * f_0)^2(y)f_0(y)\,\mathrm{d}y \leq \|K\|_1^2 \|f_0\|_2^2 \|f_0\|_\infty,$$

by Young's inequalities, and likewise for $K_{h_+}$, we have

$$\operatorname{Var}(\pi_1 L) \leq 4\|K\|_1^2 \|f_0\|_2^2 \|f_0\|_\infty =: D_1.$$

Moreover, again by Young's inequalities,

$$\|\pi_1 L\|_\infty \leq 4\|K\|_1 \|f_0\|_\infty := D_2.$$

Hence, Bernstein's inequality gives

$$\Pr(|U_n^{(1)}(\pi_1 L)| > \tau C_{n,g}/2) \leq 2\exp\left\{-\frac{n\tau^2 C_{n,g}^2/4}{2D_1 + (2/3)D_2\tau C_{n,g}/2}\right\}.$$

Since $g \geq n^{-2}(\log n)^4$, $C_{n,g} \to 0$ as $n \to \infty$ and since $g \leq h_2 = n^{-1}\ell(n)$, we have $nC_{n,g}^2 \geq d^2(g)/\ell(n)$, where we recall that $\ell(n) \to 0$, so we obtain, for any given $\tau$, that there exist $N_\tau$ such that, for all $n \geq N_\tau$,

$$\Pr(|U_n^{(1)}(\pi_1 L)| > \tau C_{n,g}/2) \leq 2\exp\{-d^2(g)/M\}. \tag{12}$$

*Second-order term:* We first state the inequality for canonical U-statistics that we are going to use (Theorem 3.4 in Houdré and Reynaud-Bouret (2003)): Let $R(x,y)$ be a symmetric function of two variables such that $ER(X,x) = 0$ for all $x$ and let

$$\Lambda_1^2 = \frac{n(n-1)}{2}ER^2,$$
$$\Lambda_2 = n\sup\{E(R(X_1,X_2)\zeta(X_1)\xi(X_2):E\zeta^2(X_1) \leq 1, E\xi^2(X_1) \leq 1\},$$
$$\Lambda_3 = \|nE_{X_1}R^2(X_1,\cdot)\|_\infty^{1/2}, \qquad \Lambda_4 = \|R\|_\infty.$$

Then, for every $\varepsilon > 0$, there exist finite non-zero numbers $\eta(\varepsilon)$, $\beta(\varepsilon)$ and $\gamma(\varepsilon)$ such that the following is true for all $u > 0$ and $n \in \mathbb{N}$:

$$\Pr\left(\frac{n(n-1)}{2}|U_n^{(2)}(R)| > 2(1+\varepsilon)^{3/2}\Lambda_1 u^{1/2} + \eta(\varepsilon)\Lambda_2 u + \beta(\varepsilon)\Lambda_3 u^{3/2} + \gamma(\varepsilon)\Lambda_4 u^2\right)$$
$$\leq 6\exp\{-u\}.$$

We apply this inequality for $R = \pi_2 L$ and $u = d^2(g)/M$ to obtain the desired bound for $\Pr(|U_n^{(2)}(\pi_2 L)| > (1-\tau)C_{n,g})$, with a small $\tau$ to be chosen below. So, we need to show that

$$2(1+\varepsilon)^{3/2}\Lambda_1 u^{1/2} + \eta(\varepsilon)\Lambda_2 u + \beta(\varepsilon)\Lambda_3 u^{3/2} + \gamma(\varepsilon)\Lambda_4 u^2 \leq (1-\tau)\frac{n(n-1)}{2}C_{n,g}$$

for the specified choice of $u$. First, since

$$\iint K_g^2(x-y)f_0(x)f_0(y)\,\mathrm{d}x\,\mathrm{d}y = \int (K_g^2 * f_0)(x)f_0(x)\,\mathrm{d}x$$
$$\leq \|f_0\|_2^2 \|K_g^2\|_1 = g^{-1}\|K\|_2^2\|f_0\|_2^2$$

and, likewise, if $g$ is replaced by $h_+ > g$, we obtain

$$\Lambda_1^2 \leq 2n(n-1)g^{-1}\|K\|_2^2\|f_0\|_2^2.$$

Taking $\varepsilon$ and $\tau$ so that $(1+\varepsilon)^{3/2} = 1.1$ and $12(1-2\tau) = 11.4$, it follows that, for all $n$ such that $\sqrt{n/(n-1)} \leq 1.1$,

$$2(1+\varepsilon)^{3/2}\Lambda_1 u^{1/2} < (1-2\tau)\frac{n(n-1)}{2}C_{n,g}.$$

For the second term,

$$|E[(K_g(X_1 - X_2))\zeta(X_1)\xi(X_2)]| \leq \|K_g * (\zeta f_0)\|_2 \|\xi f_0\|_2 \leq \|K\|_1 \|f_0\|_\infty.$$

Similarly

$$|E[E_{X_1} K_g(X_1 - X_2)\zeta(X_1)\xi(X_2)]| \leq \|K\|_1 \|f_0\|_\infty,$$
$$|EK_g(X_1 - X_2)| \leq \|K\|_1 \|f_0\|_\infty$$

and also for $K_{h^+}$. Thus,

$$E[\pi_2 L(X_1, X_2)\zeta(X_1)\xi(X_2)] \leq 8\|K\|_1 \|f_0\|_\infty$$

so that $\Lambda_2 \leq 8\|K\|_1 \|f_0\|_\infty n$. This gives that

$$\eta(\varepsilon)\Lambda_2 u = o(n^2 C_{n,g})$$

since $\sqrt{g} d(g) \to 0$ as $n \to \infty$. For the third term, we have that for every $x \in \mathbb{R}$,

$$n|E_{X_1}(\pi_2 L)^2(X_1, x)| \leq 4n[\|K\|_2^2 \|f_0\|_\infty g^{-1} + \|K\|_1^2 \|f_0\|_2^2 \|f_0\|_\infty].$$

Then,

$$\beta(\varepsilon)\Lambda_3 u^{3/2} \leq C\sqrt{n/g} d^3(g),$$

which is $o(n^2 C_{n,g})$ because $\sqrt{n}/d^2(g) \to \infty$. As for the last term, we have $\Lambda_4 = \|\pi_2 L\|_\infty \leq 4\|K\|_\infty/g$ and hence $\Lambda_4 u^2 \leq C d^4/g$, which is also $o(n^2 C_{n,g})$ because $d^3(g)$ is of the order of $(\log n)^{3/2}$, whereas $n\sqrt{g} \geq n\sqrt{h_{\min}} = (\log n)^2$. We conclude that for the specified $\tau$ and for all $n$ large enough,

$$\Pr(|U_n^{(2)}(L)| > (1 - \tau)C_{n,g}/2) \leq 6\exp\{-d^2(g)/M\}. \tag{13}$$

Inequalities (12) and (13) give

$$\Pr(|T_n(g) - ET_n(g) - (T_n(h_+) - ET_n(h_+))| > \tfrac{1}{2}\sigma(g, n) d(g)) \leq 8\exp\{-d^2(g)/M\}.$$

The lemma now follows from this bound, (11), the fact that if $g \leq h$ then

$$\exp\{-d^2(g)/M\} \leq \exp\{-d^2(h)/M\}$$

and the definition of the grid $\mathcal{H}$. □

## Acknowledgements

We would like to thank Béatrice Laurent, Alexander Tsybakov and two anonymous referees for helpful remarks.